\documentclass[12pt]{amsart}

\usepackage{epsfig}
\usepackage{amscd}
\usepackage[mathscr]{eucal}
\usepackage{amssymb}
\usepackage{amsxtra}
\usepackage{amsmath}
\usepackage[all]{xy}

\theoremstyle{plain}

\theoremstyle{definition}

\theoremstyle{remark}

\begin{document}

\bigskip

\title[An experimental study of the monotonicity of $\vert \zeta \vert$]{An experimental study of the monotonicity property of the Riemann zeta function}

\date{\today}

\author{Yochay Jerby}


%
%
% ...

\begin{abstract} In 1970, based on newly available empiric evidence, a remarkable monotonicity property for $\vert \zeta(z) \vert$ was conjectured by R. Spira. The $\zeta$-monotonicity property can be written as follows: $$ \vert \zeta (x_2 + y i ) \vert < \vert \zeta \left ( x_1 +y i \right )\vert \hspace{0.5cm} \textrm {for any } \hspace{0.25cm} x_1 < x_2 \leq 0.5 \textrm{ and } 6.29 <y. $$ In this work we present an experimental study of the monotonicity conjecture, in the course of which new properties of $\zeta(z)$ are discovered. For instance, the spectrum of semi-limits $ \lambda(z) \subset \mathbb{R}$ and the core function $C(z)$, which serves as a non-chaotic simplification of $\zeta(z)$ to the left of the critical line
\end{abstract}

\maketitle

%
%

\section{Introduction - The Riemann hypothesis and monotonicity}

\hspace{-0.6cm} In 1970, based on newly available empiric evidence,  a remarkable monotonicity property for $\vert \zeta(z) \vert$ was conjectured by R. Spira:

\bigskip

\hspace{-0.6cm} \bf The $\zeta$-monotonicity conjecture \rm (\cite{S1}): For any $y>6.29$ the function $\vert \xi \left ( x +y i \right )\vert$ is strictly-decreasing in the half-line $x<0.5$.

\bigskip

\hspace{-0.6cm} Clearly, monotonicity implies the Riemann hypothesis. In fact, the two conjectures have been shown to be equivalent, see \cite{S1,SZ,MSZ} and \cite{SC} for an analog for Riemann's $\xi$-function. Figure 1 illustrates the $\zeta$-monotonicty property in the domain $0<y< 10^4$:

\begin{figure}[h!]
\centering
\includegraphics[width=0.75 \linewidth]{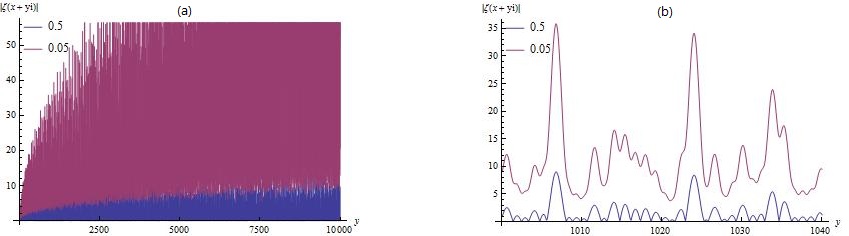}
\caption{Graph of $\vert \zeta(0.5+yi) \vert $ and $\vert \zeta (0.05+yi) \vert $ for $0 \leq y \leq 10^4$ (a) and $10^3 \leq y \leq 10^3+40$ (b), showing the property $ \vert \zeta(0.5+yi) \vert < \vert \zeta(0.05+yi) \vert$. }
\label{fig:3}
\end{figure}

\hspace{-0.6cm} This work is devoted to a further, modern, experimental study of the $\zeta$-monotonicity conjecture and related properties. In particular, this study leads to the discovery of various other new fundamental properties of zeta. For instance: the \emph{spectrum of semi-limits} of the partial sums, which we view as governing the chaotic part of zeta, and to the definition of the \emph{core}, $C(z)$, which serves as a non-chaotic simplification of $\zeta(z)$, to the left of the critical line.

\section{A few effective remarks on the monotonicity property}
\label{s:left}

\hspace{-0.6cm} Let us consider the function
\begin{equation} \label{eq:4}
\eta(y,t) := e^t \left ( \vert \zeta(0.5 \cdot (1-e^{-t}) +yi ) \vert - \vert \zeta(0.5+yi) \vert \right ).
\end{equation}
The $\zeta$-monotonicity implies that $\eta(y,t)$ is positive for $y \geq 6.29 $ and $ t \geq 0$ or, equivalently, that the function $log ( \eta(y,t))$ is well-defined. Figure 2 shows, for instance, the values of $log ( \eta(y,t)) $ for $t=0$ and $t=10$ in the domain $6.29 \leq y \leq 2 \cdot 10^3 $:
\begin{figure}[h!]
\centering
\includegraphics[scale=0.35]{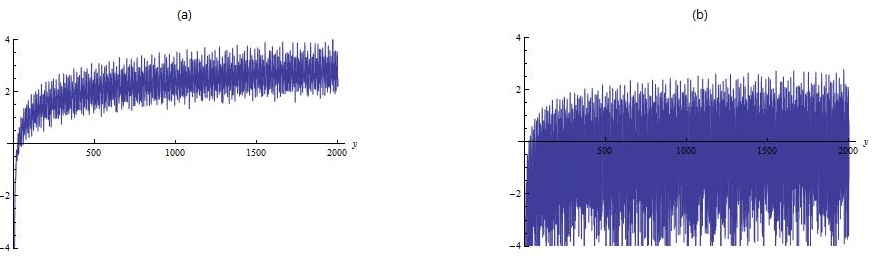}
\caption{Graph of $log (\eta(y,0))$ (a) and $log( \eta(y,10))$ (b) for $6.29 \leq y \leq 2 \cdot 10^3$.}
\end{figure}

\hspace{-0.6cm} Note that the function $log(\eta(y,0))$ seems to be not only well-defined but, also, to fluctuate, at its core, around some strictly increasing function. However, as $t$ increases, it becomes less straight-forward to discern that the $log(\eta(y,t))$ is well define. However, it turns that the function $log (\eta(y,t))$ also admits the following remarkable property: it could be bounded from below in terms of $log(\vert \zeta(0.5+yi) \vert)$ itself, in the region $y>6.29$ (compare, for instance, Fig. 3)!
\begin{figure}[h!]
\centering
\includegraphics[scale=0.35]{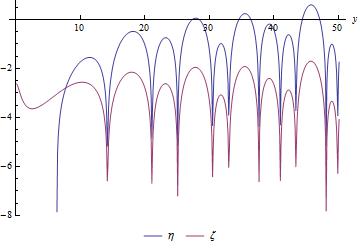}
\caption{Graph of $log (\eta(y,10))$ and $log (\vert \zeta(0.5+yi) \vert ) -3$ over $0 \leq y \leq 50$. \label{overflow}}
\end{figure}

\hspace{-0.6cm} This leads us to consider, instead of $\eta(y,t)$, the following function
\begin{equation} \label{eq:5}
\widetilde{\eta} (y,t):= e^{t} \cdot \left ( \frac{ \vert \zeta (0.5(1-e^{-t})+yi) \vert }{ \vert \zeta(0.5 +yi) \vert } -1 \right ),
\end{equation} which is defined for all $(y,t)$ such that $\zeta(0.5+yi) \neq 0 $ and positive exactly when $\eta(y,t)$ is. The advantage of $\widetilde{\eta} (y,t)$ over $\eta(y,t)$ is that, contrary to $\eta(y,t)$, the function $\widetilde{\eta} (y,t)$ appears to be not only strictly positive, but actually seems to be bounded from below by a rather well behaved, \emph{smooth, non-chaotic, increasing} function $\widetilde{X}(y,t)$, for any given $t$ (as Fig. 4 illustrates).

\begin{figure}[h!]
\centering
\includegraphics[scale=0.45]{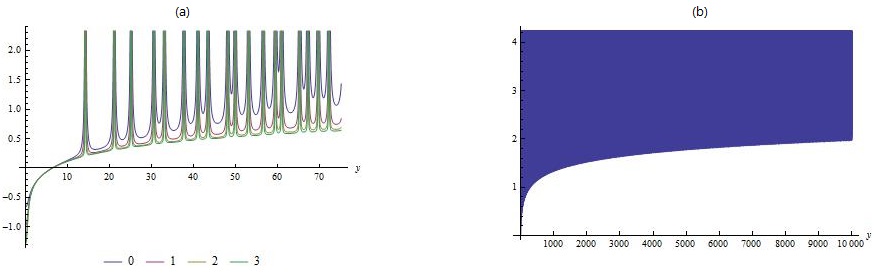}
\caption{Graph of $\widetilde{\eta} (y,t)$ over $0<y<75$ and $t=0,1,2,3$ (a) and $\widetilde{\eta} (y,3)$ over $0 \leq y \leq 10^4$ (b).}
\end{figure}

\hspace{-0.6cm} Recall that the functional equation of $\zeta(z)$ is given by $ \zeta(z) = \chi (z) \cdot \zeta(1-z)$ where
$
\chi(z) := 2^z \pi^{z-1} sin \left ( \frac{\pi z }{2} \right ) \Gamma(1-z)$. It turns that a rather good first order approximation of $\widetilde{X}(y,t)$ could be given in terms of the following function: \begin{equation} \label{eq:7}
X(y,t):= e^{t} \cdot \left ( \frac{ \vert \chi (0.5(1-e^{-t})+yi) \vert }{ \vert \chi(0.5 +yi) \vert } -1 \right ) = e^t \cdot \left ( \vert \chi (0.5(1-e^{-t})+yi) \vert -1 \right ).
\end{equation} Figure 5, for instance, shows a graph of $log(\widetilde{\eta} (y,10) ) $ (blue) and $log(X(y,10))-0.75$ (purple) for $0<y<10^4$:
\begin{figure}[h!]
\centering
\includegraphics[scale=0.45]{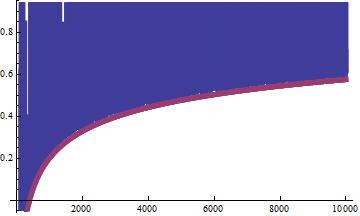}
\caption{Graph of $log(\widetilde{\eta}(y,10))$ and $log(X(y,10))-0.75$ over $0<y<10^4$.}
\end{figure}

\hspace{-0.6cm} Our aim is to explain why the increasing non-chaotic $\widetilde{X}(y,t)$ such that $\widetilde{X}(y,t) < \widetilde{\eta} (y,t)$ should countinue to exist for any $y \geq 6.29$ and $t \geq 0$. In order to do so, we need to introduce the spectrum of semi-limits and core function in the next section 3.

\section{The spectrum of semi-limits and the core function}
\label{s:left}

\hspace{-0.6cm} Recall that in the critical strip $0 \leq Re(z) \leq 1$ zeta is given by $\zeta(z)=\frac{1}{1-2^{1-z}} \sum_{k=1}^{\infty} \frac{(-1)^{k+1}}{k^z}$.
For $n \in \mathbb{N}$ consider the partial sums \begin{equation} \label{eq:10} S_n(z):=\frac{1}{1-2^{1-z}} \sum_{k=1}^{n} \frac{(-1)^{k+1}}{k^z}.
\end{equation}

\hspace{-0.6cm} The starting point of this section is the observation of a few special properties of the series $S_n(z)$. For instance, a typical example of the behavior of $S_n(z)$ is presented in Fig. 6:

\begin{figure}[h!]
\centering
\includegraphics[scale=0.4]{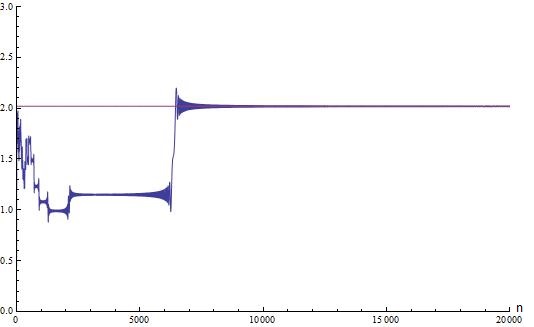}
\caption{Values of $ \vert S_n(0.5 \cdot (1-e^{-2}) +2 \cdot 10^4 \cdot i) \vert$ for $n=0,...,2 \cdot 10^4$. }
\end{figure}

\hspace{-0.6cm} As one can see, the series $S_n(z)$ actually fluctuates around various other values for a while before "starting to approximate" $\zeta(z)$ (purple) and that the "surge towards $\zeta(z)$" is made around the "critical" stage $n \approx Im(z)/3$. In fact, the bigger $Im(z)$ becomes, the interval $[0,Im(z)/3]$ becomes divided into more and more sub-segments over which $S_n(z)$ fluctuates around a certain fixed \emph{semi-limit}, and the transition between two such \emph{semi-limits} is done by steep surges (as in the picture). We refer to the collection of these semi-limits $\lambda(z) \subset \mathbb{R}$ as the \emph{spectrum of the value $z$}.

\hspace{-0.6cm} In particular, in view of the above, for $\alpha \in [0,1]$, we define the \emph{$\alpha$-truncation of zeta}:
\begin{equation} \label{eq:11} \zeta_{\alpha}(z):=\frac{1}{1-2^{1-z}} \sum_{k=1+ [(1-\alpha) \cdot Im(z)]}^{\infty} \frac{(-1)^{k+1}}{k^z},
\end{equation} where $[y] \in \mathbb{N}$ stands for the \emph{integral value} of the real number $y$. An example is presented in Fig. 7:
\begin{figure}[h!]
\centering
\includegraphics[scale=0.45]{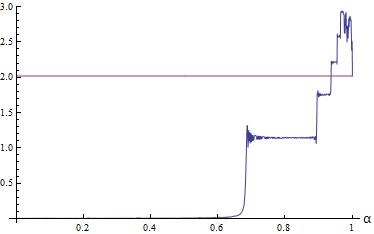}
\caption{Values of $ \vert \zeta_{\alpha}(0.5 \cdot (1-e^{-2})+2 \cdot 10^4 \cdot i) \vert$ for $0 \leq \alpha \leq 1$.}
\end{figure}

\hspace{-0.6cm} It is very interesting to study the properties of the \emph{spectrum} of \emph{semi-limits} which the $\alpha$-truncation $\zeta_{\alpha}(z)$ fluctuates around, in general (see remark 3.1 below). However, of special interest for us is the first of them (the last being the limit $\zeta(z)$ itself). In particular, truncating from $\zeta(z)$ all the semi-limits except for the first one leads to the following definition:

\bigskip

\hspace{-0.6cm} \bf Definition 3.1: \rm $C(z):=\vert \zeta_{0.8}(z) \vert $ is the core function of $\zeta(z)$.

\bigskip

\hspace{-0.6cm} Let us note that the value $a=0.8$ is simply taken to represent the first spectrum value, which occurs for $\zeta_{\alpha}(z)$ around $\alpha \approx 2/3$. We call $C(z)$ the core of $\zeta(z)$ as, left to the critical line, the core turns to serve as a non-chaotic simplification of $\vert \zeta(z) \vert$, as illustrated in Fig. 8:
\begin{figure}[h!]
\centering
\includegraphics[scale=0.35]{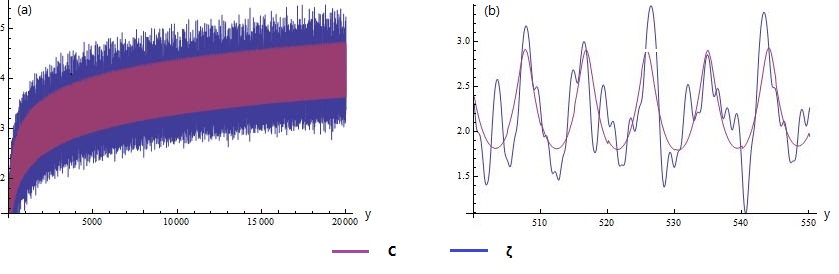}
\caption{Graph of $log(C(yi))$ and $log \vert \zeta(yi) \vert$ for $0 \leq y \leq 2 \cdot 10^4$ (a) and $500 \leq y \leq 550$ (b). }
\end{figure}

\hspace{-0.6cm} Let us note that, as by truncating the semi-limits in the spectrum $\lambda(z)$ we obtained the non-chaotic
core function $C(z)$, we view $\lambda(z)$ as encoding the chaotic, random, features of zeta. Moroever, as mentioned the number of elements in $N(z) = \vert \lambda(z) \vert$ grows as $Im(z) \rightarrow \infty$ (see further discussion in section 4).

\hspace{-0.6cm} The importance of the core $C(z)$ to the study of the $\zeta$-monotonicity property is that it can be viewed as the part of $\zeta(z)$ that is efficiently approximated in terms of $\vert \chi(z) \vert$. Figure 9 illustrates this for $Re(z)=0$:

\begin{figure}[h!]
\centering
\includegraphics[scale=0.4]{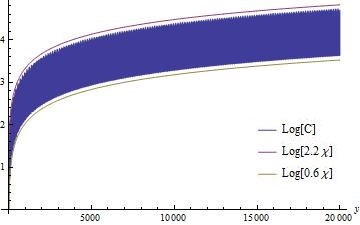}
\caption{Graph of $log(C(yi))$ and $ log(2.2 \vert \chi (yi) \vert) ,log(0.6 \chi(yi) \vert )$, approximating $log(C(yi))$ from above and below, for $0 \leq y \leq 2 \cdot 10^4$. }
\end{figure}
\hspace{-0.6cm} In particular, the "mysterious difference" between $\widetilde{X}(y,t)$ and $X(y,t)$ can be viewed as a result of the contribution of the re-addition of the chaotic elements of the spectrum $\lambda(z)$. Moreover, recall that the core, $C(z)= \vert \zeta_{0.8}(z) \vert$, by definition, conceptually represents $80 \%$ of the relevant elements of the series defining $\vert \zeta(z) \vert$. Hence, we turn to discuss the addition the remaining, chaotic, $20 \%$, to $\zeta(z)$ and, mainly, to $\widetilde{\eta} (y,t)$. For $ a \in [0,1]$ set:
\begin{equation} \label{eq:12}
C_a(z) :=\left \vert \zeta_{0.8}(z) + \frac{a}{(1-2^{1-z} ) } \cdot \sum_{k=1}^{[0.2Im(z)]} \frac{ (-1)^{k+1}}{k^z} \right \vert,
\end{equation} Interpolating between the core, $C(z)=C_0(z)$, and zeta itself, $\vert \zeta (z) \vert = C_1(z)$. In view of section 2 set:
\begin{equation} \label{eq:13}
\widetilde{\eta} _a(y,t) := e^t \cdot \left ( \frac{C_a(0.5(1-e^{-t})+yi)}{C_a(0.5+yi)}-1 \right ).
\end{equation}
Figure 10 shows the remarkably structured way $\widetilde{\eta} _a(y,t)$ transitions from $\widetilde{\eta} _0(y,t)$ (blue) to $\widetilde{\eta} (y,t)=\widetilde{\eta} _1(y,t)$ (purple) (contrary to the chaotic transition of $C(z)$ to $\vert \zeta(z) \vert$):

\begin{figure}[h!]
\centering
\includegraphics[scale=0.4]{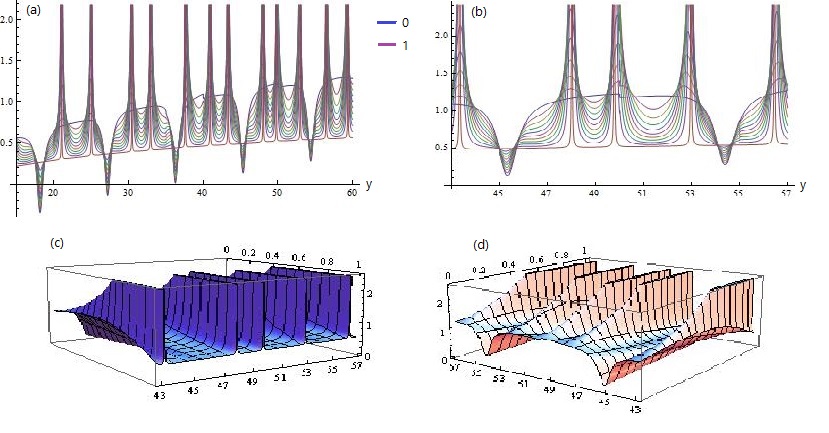}
\caption{$\widetilde{\eta}_{k/10} (y,5)$ with $k = 0,...,10$ for $15 \leq y \leq 60 $ (a) and $43 \leq y \leq 57$ (b) and of $\widetilde{\eta} _a (y,5)$ over $(y,a) \in [43,57] \times [0,1] $ front (c) rear (d).}
\end{figure}

\hspace{-0.6cm} The key feature is that, on the one hand, the structured transition from $\widetilde{\eta} _0(y,t)$ to $\widetilde{\eta} _1(y,t)$, described in Fig. 10, is independent of $y$. However, on the other hand, by definition, the number of elements discerning between $\widetilde{\eta} _0(y,t)$ and $\widetilde{\eta} _1(y,t)$ is given by $[0.2 y]$. Hence, conceptually, the pattern in Fig. 10 can be explained in terms of "induction on $[0.2y]$".

\hspace{-0.6cm} Finally, let us note the following remark regarding the poles: Set
\begin{equation} \label{eq:14} \begin{array}{ccc} \eta(y,t) : =e^t ( \vert \zeta(0.5(1-e^{-t})+yi) \vert - \vert \zeta(0.5 +yi) \vert )
& ; & \theta(y):= \vert \zeta (0.5 +yi) \vert. \end{array}
\end{equation}
In fact, it seems possible to locally \emph{resolve} the poles altogether by replacing $ \theta(y) $ with a smooth non-vanishing function $\widetilde{\theta} (y) \neq 0$, coinciding with $\theta(y)$ away from small neighborhoods of the zeros, and keeping the property $\widetilde{X} (y,t) < \widetilde{\widetilde{\eta} }(y,t):= \eta(y,t) / \widetilde{\theta}(y) $. In order to understand how the local correction should occur let us consider $\widetilde{y}_1 \approx 14.1347$, the first zero on the critical strip. Figure 11 shows the behavior of $\eta(y,t)$ and $\theta(y)$ in a small neighborhood of $\widetilde{y}_1$:

\begin{figure}[h!]
\centering
\includegraphics[scale=0.4]{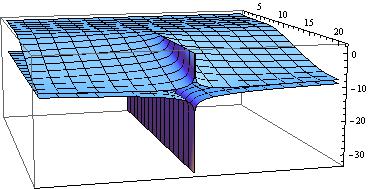}
\caption{Graph of $log(\theta(\widetilde{y}_1+ \epsilon)), log(\eta(\widetilde{y}_1+ \epsilon,t))$ for $\vert \epsilon \vert \leq 0.0001$ and $0 \leq t \leq 20 $.}
\end{figure}

\hspace{-0.6cm} In particular, we can take $\widetilde{\theta}(y) = max ( \theta(y) , e^{-20})$ as the required correction in the considered neighborhood of $\widetilde{y}_1$. Let $\widetilde{z}_k = 0.5+\widetilde{y}_k$ be the $k$-th zero of zeta on the critical strip. Empirical verification shows that the typical local behavior of $log(\eta(y,t))$ and $ log( \theta(y))$ in a neighborhood of $\widetilde{y}_k$ for any $k$ is, in fact, similar to that presented in Fig. 11 for $\widetilde{y}_1$.

\bigskip

\hspace{-0.6cm} \bf Remark 3.2 \rm (The right-hand-side): The function $log \vert \zeta( x +yi) \vert $ is known to be unbounded, as a function of $y$, for $0.5 \leq x \leq 1$, see, for instance, Theorem 11.9 of \cite{T}. In view of this it is interesting to point out Fig. 12:

\begin{figure}[h!]
\centering
\includegraphics[scale=0.4]{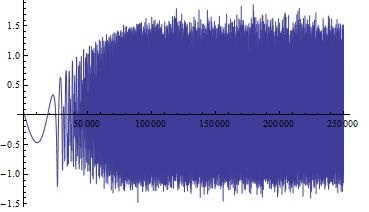}
\caption{Graph of $log \vert \zeta(0.95+e^{0.0001t} i) \vert$ over $t=0,...,250000$.}
\end{figure}

\hspace{-0.6cm} First, as one can see from figure 12, even though $log \vert \zeta(0.95 +yi) \vert$ is guaranteed to be unbounded it should nevertheless have an extremely slow rate of growth. However, in view of Fig. 12 it also becomes interesting to ask the following more refined global question:

\bigskip

\hspace{-0.6cm} \bf Question: \rm For $0.5 \leq x \leq 1$ and $ s \in \mathbb{R}$ let $Y(s,x)$ be the minimal value of $0 \leq Y$ for which $log \vert \zeta(x +Yi) \vert = s$. What can be said about the function $Y(s,x)$?

\bigskip

\hspace{-0.6cm} For instance, Fig. 12 shows that $e^{25}<Y(\pm 2,0.95)$. The above question has, of course, direct bearing on zeros of zeta. As mentioned, the results of \cite{T} and later, the, much more general, Voronin's universality theorem \cite{V} imply that $log \vert \zeta(x+yi) \vert$ is unbounded in $0.5 \leq x \leq 1$. However, these classical results are non-quantitative, in the sense that they do not give quantitative information on $Y(s,x)$ beyond guaranteeing that it is unbounded. In particular, this does not exclude the existence of, so called, "ghost zeros", that is, an infinite amount of extremely tiny non-zero values which (a) can hardly be discerned from a real zero (b) appear practically anywhere to the right of the critical line. However, Fig. 15 shows that this is not exactly the case. Indeed, at least for $0<y<e^{25}$ the size of $\vert \zeta (0.95+yi) \vert $ is globally bounded by $e^{-2}$ and, as $y$ grows, it is natural to suggest that it would be possible to extend this bound by a (very slowly) decreasing function of $y$. It is important to note in this context results of Garunkštis on effective versions of Voronin's universality theorem, specifically corollary 2 of \cite{G}, which also seem to suggest very slow asymptotics for $Y(s,x)$. Moreover, in the context of this work, it is interesting to note that the question of the description of $Y(s,x)$ could be viewed as the right hand side analog of the description of the "mysterious difference" between $\widetilde{X}(y,t)$ and $X(y,t)$.

\section{Concluding remarks}
\label{s:intro}

\hspace{-0.6cm} In this work we conducted an experimental study of the $\zeta$-monotonicity conjecture \cite{BB,S1,SZ}, which is an equivalent reformulation of the Riemann hypothesis, see \cite{MSZ,S1,SC}. This led us to discover the spectrum of semi-limits $\lambda(z) \subset \mathbb{R}$ (which we view as dominating the chaotic features of zeta) and the existence of the core function $C(z)$ which we view as a non-chaotic simplification of $\vert \zeta(z) \vert$ (to the left of the critical line), obtained by truncating the semi-limits in $\lambda(z)$ aside from the first one. As mentioned, for a given $z \in \mathbb{C}$, the spectrum is a collection of random-like values such that $N(z) = \vert \lambda (z) \vert \rightarrow \infty $ when $Im(z) \rightarrow \infty$. One of the fascinating aspects in the modern study of zeta is the discovery of various relations to quantum chaos, see \cite{Bog} and references therein. Specifically the existence of conjectural relations between statistical properties of the zeros of zeta and statistical properties of $\lambda(M)$, the spectrum of eigenvalues of random $N \times N$-matrices in GUE (Gaussian Unitary Ensemble), such that $N \rightarrow \infty$. Even though, there is vast empirical evidence to back up the the various quantum chaos conjectures, a conceptual explanation to the observed relation between zeros of zeta and eigenvalues of random matrices, is still largely missing. In view of this, it is interesting to ask, weather the monotonicity reformulation can be extended to relate zeros of zeta, the spectrum of semi-limitis $\lambda(z)$ and spectrum of eigenvalues of random matrices $\lambda(M)$?

\end{document}